\documentclass[a4paper,12pt]{amsart}
\usepackage{ifthen}
\usepackage{graphicx}
\usepackage{mathrsfs}
\usepackage{color}
\nonstopmode \numberwithin{equation}{section}
\makeatletter
\@namedef{subjclassname@2020}{\textup{2020} Mathematics Subject Classification}
\makeatother

\newtheorem{thm}{Theorem}[section]
\newtheorem{cor}[thm]{Corollary}
\newtheorem{lem}[thm]{Lemma}
\newtheorem{prop}[thm]{Proposition}
\newtheorem{rem}[thm]{Remark}

\newtheorem{prob}[thm]{Problem}

\theoremstyle{definition}
\newtheorem{defn}{Definition}

\newenvironment{pf}[1][]{%
 \vskip 3mm
 \noindent
 \ifthenelse{\equal{#1}{}}%
  {{\slshape Proof. }}%
  {{\slshape #1.} }%
 }%
{\qed\bigskip}

\newcounter{alphabet}

\newenvironment{Thm}[1][]{\refstepcounter{alphabet}%
\bigskip%
\noindent%
{\bf Theorem \Alph{alphabet}}%
\ifthenelse{\equal{#1}{}}{}{ (#1)}%
{\bf .} \itshape}{\vskip 8pt}

\newcommand{\C}{{\mathbb C}}
\newcommand{\D}{{\mathbb D}}

\newcommand{\PP}{{\mathcal P}}
\newcommand{\R}{{\mathbb R}}

\newcommand{\T}{{\mathcal T}}

\renewcommand{\Im}{{\,\operatorname{Im}\,}}

\renewcommand{\Re}{{\,\operatorname{Re}\,}}

\newcommand{\Gauss}{{\null_2F_1}}

\renewcommand{\arg}{\,{\operatorname{arg}\,}}

\newcommand{\aand}{{\quad\text{and}\quad}}

\newcounter{minutes}
\divide\time by 60
\newcounter{hours}
\multiply\time by 60 \addtocounter{minutes}{-\time}

\begin{document}
\bibliographystyle{amsplain}
\title[Universal convexity and range problems]%
{Universal convexity and range problems of shifted hypergeometric functions
}

\begin{center}
{\tiny \texttt{FILE:~\jobname .tex,
        printed: \number\year-\number\month-\number\day,
        \thehours.\ifnum\theminutes<10{0}\fi\theminutes}
}
\end{center}
\author[T. Sugawa]{Toshiyuki Sugawa}
\address{Graduate School of Information Sciences,
Tohoku University,
Aoba-ku, Sendai 980-8579, Japan}
\email{sugawa@math.is.tohoku.ac.jp}
\author[L.-M.~Wang]{Li-Mei Wang}
\address{School of Statistics,
University of International Business and Economics, No.~10, Huixin
Dongjie, Chaoyang District, Beijing 100029, P.R. China}
\email{wangmabel@163.com}
\author[C.F.~Wu]{Chengfa Wu}
\address{Institute for Advanced Study, Shenzhen University, Shenzhen 518060,
P. R. China\\
School of Mathematical Sciences, Shenzhen University, Shenzhen 518060,
P. R. China}
\email{cfwu@szu.edu.cn}


\begin{abstract}
In the present paper, we study the shifted
hypergeometric function $f(z)=z\Gauss(a,b;c;z)$ for real parameters
with $0<a\le b\le c$ and its variant $g(z)=z\Gauss(a,b;c;z^2).$
Our first purpose is to solve the range problems for $f$ and $g$ posed by
Ponnusamy and Vuorinen in their 2001 paper.
Ruscheweyh, Salinas and Sugawa developed in their 2009 paper
the theory of universal prestarlike functions
on the slit domain $\C\setminus[1,+\infty)$ and showed universal
starlikeness of $f$ under some assumptions on the parameters.
However, there has been no systematic study of universal convexity of
the shifted hypergeometric functions except for the case $b=1.$
Our second purpose is to show universal convexity of $f$
under certain conditions on the parameters.
\end{abstract}

\thanks{The second author is supported by a grant of University of International Business and Economics (No. 78210418) and National Natural Science Foundation of China (No.
11901086) .
The third author was supported by the National
Natural Science Foundation of China (grant numbers 11701382 and 11971288); and
Guangdong Basic and Applied Basic Research Foundation, China (grant number
2021A1515010054).
}
\maketitle

\section{Introduction and main results}

For real numbers $a,~b,~c$ with $c\not=0,-1,-2,\ldots$,
the Gaussian hypergeometric series
$$
\Gauss(a,b;c;z)=\sum_{n=0}^{\infty}\frac{(a)_n(b)_n}{(c)_nn!}z^n
$$
defines a function analytic on $\D:=\{z\in\C: |z|<1\}$.
Here the \emph{Pochhammer symbol} $(a)_n$
denotes the ascending (or shifted) factorial function
for $n\in\mathbb{N}:=\{0,1,2,\dots\};$ that is,
$(a)_0=1, ~(a)_{n}=a(a+1)\cdots(a+n-1)=\Gamma(a+n)/\Gamma(a)$
for $n=1,2,\dots.$
Hypergeometric functions can be
analytically continued as a solution to the hypergeometric differential equation
along any path in the complex plane that avoids the
regular singularities $0$, $1$ and $\infty$ of the equation.
In particular, we can regard $\Gauss(a,b;c;z)$ as an analytic function
on the slit plane $\Lambda=\{z\in\C: |\arg(1-z)|<\pi\}=
\C\setminus [1,+\infty)$.
For basic properties of hypergeometric functions we refer the readers to
\cite{AbramowitzStegun:1965}, \cite{OLBC:2010} and \cite{Wall:anal}.
In the present paper, we will study geometric properties of the
\emph{shifted} hypergeometric function $z\Gauss(a,b;c;z).$
Since $\Gauss(a,b;c;z)=\Gauss(b,a;c;z),$ we may assume $a\le b$ without
loss of generality.

The behavior of hypergeometric functions $\Gauss(a,b;c;z)$
near $z=1$ can be classified into three cases according to the
sign $a+b-c$ (see Lemma \ref{lem:1-} below for details).
In particular, hypergeometric functions are unbounded for the case $c\leq a+b$.
If $c<a+b$, Ponnusamy and Vuorinen \cite[p.~351]{PonVuor:2001}
posed the following question (see also \cite{HPV:2010}):

\begin{prob}
\label{prob}
Suppose $a, ~b,~ c>0$ with $c<a+b$.
Do there exist $\delta_1, ~\delta_2>0$ such that
for $a \in\left(0, \delta_1\right)$ and $b \in\left(0, \delta_2\right)$
the function $z\Gauss(a, b ; c;  z)~ (z\Gauss\left(a, b ; c ; z^2\right)$ respectively$)$
has the property that the image domain of $\D$ is completely contained in a sector type domain
where the ``angle" depends on $a+b-c$?
\end{prob}

This paper is motivated partly by this problem.
Note that the balanced case $a+b=c$ was treated in \cite{Wang2022}.
First, we need to formulate sector type domains.
For $0<\alpha<1,$
we denote by $S(\alpha)$ the sector $\{z: |\arg z|<\pi\alpha/2\}$ of opening angle $\pi\alpha$
and set $S^*(\alpha)=\{z: z\in S(\alpha) ~\text{or}~ -z\in S(\alpha)\}
=S(\alpha)\cup[-S(\alpha)].$
Moreover, for $\varepsilon>0,$ we denote by $S_\varepsilon(\alpha)$ and $S_\varepsilon^*(\alpha)$
the $\varepsilon$-neighbourhoods of $S(\alpha)$ and $S^*(\alpha),$ respectively.
That is, $S_\varepsilon(\alpha)=\{z\in\C: |z-z'|<\varepsilon ~\text{for some}~ z'\in S(\alpha)\}$
and $S_\varepsilon^*(\alpha)$ is defined likewise.

We answer to the above problem in the following form.

\begin{thm}\label{thm:1}
Let $a, b, c$ be positive numbers with $a\le b$ and $\delta=a+b-c\in(0,1)$ and
$$
\varepsilon=\left|\frac{\Gamma(c)\Gamma(-\delta)}{\Gamma(c-a)\Gamma(c-b)}\right|
+2^{1-\delta}\cdot\frac{\Gamma(c)\Gamma(\delta)}{\Gamma(a)\Gamma(b)}.
$$
Then the following hold.
\begin{enumerate}
\item[(I)]
The image of the unit disk $\D$ under the mapping
$z\Gauss(a, b ; c;  z)$ is contained in $S_\varepsilon(\delta).$
\item[(II)]
The image of the unit disk $\D$ under the mapping
$z\Gauss(a, b ; c;  z^2)$ is contained in $S_\varepsilon^*(\delta).$
\end{enumerate}
\end{thm}

When either $c=a$ or $c=b,$ we understand $\varepsilon=2^{1-\delta}$ in the above.
Note that the above theorem does not tell univalence of the mappings.
Indeed, Ponnusamy and Vuorinen \cite[Thm 6.1]{PonVuor:2001} observed that
$f(z)=z\Gauss(a, b ; c;  z)$ is not univalent for $a,b\in(0,2)$ and $0<c<ab/2.$
However, if $f(z)$ is univalent and has a nice geometric property, we may make a
stronger assertion.

For a function $f$ analytic on $\D$ and
normalized by $f(0)=f'(0)-1=0$,
the \emph{order of convexity} of $f$ is defined by
 $$
\kappa=\kappa(f):=\inf_{z\in \D}\Re \left(1+\frac{zf''(z)}{f'(z)}\right)\in[-\infty,1].
$$

It is known that $f$ is \emph{convex}, i.e. $f$ is univalent on $\D$ and $f(\D)$ is a convex domain
if and only if $\kappa(f)\geq 0$.
It is also known that if $\kappa(f)\geq -1/2$
then $f$ is univalent on $\D$ and $f(\D)$ is convex in
(at least) one direction, see \cite{Umezawa:1952}
and \cite[p.73]{Rus:conv}.
Starlikeness of shifted hypergeometric functions was thoroughly studied by K\"ustner
\cite{Kustner:2002, Kustner:2007}.
However, convexity of shifted hypergeometric functions is seemingly difficult to establish
except for the case when $b=1,$ see \cite[pp.~1368--1369]{Kustner:2007}.
Recently, one of the authors obtained the following result concerning the order of convexity.

\begin{Thm}$($\cite[Thm.1.1]{Wang2021}$)$
For real parameters $a, b$
and $c$ satisfying $0<a<1<b \leq c<\min\{a+b,1+a+b-ab\}$,
the order of convexity of the function $z\Gauss(a, b; c; z)$
is
$$
\kappa
=\frac{c^2-a^2-b^2+3(a+b-c)-2}{2(a+b-c)}.
$$
\end{Thm}

In particular, $z\Gauss(a, b; c; z)$ is convex when $c^2-a^2-b^2+3(a+b-c)\ge 2.$
For convex functions, we can refine Theorem \ref{thm:1} in the following form.

\begin{thm}\label{thm:sector}
For positive numbers $a,~b$ and $c$ with $\delta=a+b-c\in(0,1),$
suppose that the shifted hypergeometric function
$f(z)=z\Gauss(a, b ; c;  z)$ is convex on $\D.$
Then the image $f(\D)$ is contained in the sector $S=\{z: |\arg(z-B)|<\pi\delta/2\}$
and the boundary of $S$ is contained in the union of asymptotic lines
of the boundary of $f(\D),$ where $B=\Gamma(c)\Gamma(-\delta)/[\Gamma(c-a)\Gamma(c-b)].$
Moreover, the choice of the sector $S$ is optimal.
\end{thm}

Since the origin is contained in $f(\D),$ the inequality $B<0$ must hold
under the assumption of Theorem \ref{thm:sector}.
Note also that $S\subset S_{-B}(\delta).$

Ruscheweyh, Salinas and Sugawa in \cite{RusSalSu} extended the
study of geometric properties of analytic functions
on the unit disk to other disks and half planes containing
the origin. They introduced the class of
universally prestarlike functions of order $\alpha\leq 1$ in the domain $\Lambda$
(containing universally convex and universally starlike functions as the
special cases when $\alpha=0$ and $\alpha=1/2,$ respectively).
A general condition of a shifted hypergeometric function
to be universally starlike is given in \cite{RusSalSu}.
However, a similar condition for universal convexity is not known
so far except for the case $b=1$ (see Theorem C below)
to the knowledge of the authors.
(In \cite{RusSalSu}, universal convexity is shown for functions of the form
$(c/ab)[\Gauss(a,b;c;z)-1].$)
In the present paper,
we investigate the universal convexity of the shifted hypergeometric functions.

\begin{defn}$($\cite[Def.1.4]{RusSalSu}$)$\label{uc}
A function holomorphic in $\Lambda=\C\setminus[1,+\infty)$ is called \emph{universally convex}
if it maps every half-plane or disk in $\Lambda$ containing the origin
onto a convex domain.
\end{defn}

We denote by $\PP$ the class of analytic functions $P$ with $P(0)=1$ and
$\Re P>1/2$ on $\D.$
The following formula is known as the Herglotz representation:
$$
P(z)=\int_{\partial\D} \frac{d\mu(\zeta)}{1-\bar\zeta z},\quad z\in\D,
$$
for a (Borel) probability measure $\mu$ on the unit circle $\partial\D.$
(In what follows, we will surpress ``Borel" because we will consider only Borel measures here.)
Similarly, we consider the class $\T$ of those analytic functions $F$ with $F(0)=1$
which admit the representation
\begin{equation}\label{eq:F}
F(z)=\int_0^1 \frac{d\mu(t)}{1-tz}, \quad z\in\Lambda,
\end{equation}
for a probability measure $\mu$ on the interval $[0,1].$
Such a function $F$ is expanded in the form of power series $F(z)=c_0+c_1z+c_2z^2+\dots$ on
$|z|<1$ for a totally monotone (or completely monotone) sequence $\{c_n\}$ with $c_0=1,$ that is,
$\Delta^kc_n\ge 0$ for all $k,n\ge0,$
where $\Delta^0c_n=c_n$ and $\Delta^{k+1}c_n=\Delta^kc_n-\Delta^kc_{n+1}$
for $k,n\ge 0.$
Conversely, for a totally monotone sequence $\{c_n\}$ with $c_0=1,$ the function
$F(z)=\sum c_nz^n$ is known to be in the class $\T$ (see \cite{Wall:anal} for details).

An analytic characterization of universally convex functions is given in terms of
the class $\T$ as follows.

\begin{Thm}$($\cite[Thm.1.6]{RusSalSu}$)$\label{Thm:UnivConvex}
An analytic function $f$ on $\Lambda$ with $f(0)=0$ and $f'(0)=1$ is universally convex
if and only if $1+zf''(z)/(2f'(z))$ belongs to the class $\T.$
\end{Thm}

We note that a normalized analytic function $f$ on $\D$ is convex
if and only if $1+zf''(z)/(2f'(z))$ belongs to the class $\PP.$
Since $\T\subset\PP,$ we see that a universally convex function is always convex.

Let $f(z)=z+a_2z^2+a_3z^3+\dots$ be a universally convex function.
Then, by a formal computation, we have for $1+zf''(z)/(2f'(z))=1+c_1z+c_2z^2+\dots$
$$
c_1=a_2,\quad
c_2=3a_3-2a_2^2, \aand
c_3=6a_4-9a_2a_3+4a_2^3.
$$
Since $\{c_n\}$ with $c_0=1$ is totally monotone, we have
countably many inequalities such as
$$
0\le c_1=a_2\le 1 \aand
0\le c_1-c_2=a_2-3a_3+2a_2^2\le 1-c_1=1-a_2.
$$
For the shifted hypergeometric function $f(z)=z\Gauss(a,b;c;z)
=z+\alpha z^2+\alpha\beta z^3+\dots$ with $\alpha=ab/c$
and $\beta=(a+1)(b+1)/[2(c+1)],$ we have the necessary conditions
for universal convexity of $f$ such as
$$
0\le \alpha\le 1, \quad 3\beta\le 1+2\alpha
 \aand \alpha(2+2\alpha-3\beta)\le 1.
$$
It is, of course, not practical to check all the inequalities $\Delta^kc_n\ge0$
to obtain a sufficient condition for universal convexity of $f.$

Since $zf'(z)=z\Gauss(a,2;c;z)$ for $f(z)=z\Gauss(a,1;c;z),$ we have
the following result immediately from
Ruscheweyh et al.~\cite[Thm.1.8]{RusSalSu}.

\begin{Thm}
Let $a$ and $c$ be positive numbers with $a\le 1$ and $2\le c.$
Then the shifted hypergeometric function $f(z)=z\Gauss(a,1;c;z)$
is universally convex.
\end{Thm}

Since the theorem does not appear in this form in the literature, we will include
a brief account of the proof after Remark \ref{rem:T} below.
In the present paper, the authors extend it in the following way.

\begin{thm}\label{thm:convex}
Let $a,~b$ and $c$ be positive numbers with $a\le b\le c.$
If one of the conditions $(I)$ and $(II)$ below is satisfied,
then the shifted hypergeometric function $f(z)=z\Gauss(a, b; c; z)$ is
universally convex on the slit plane $\Lambda=\C\setminus [1,+\infty)$:
\begin{enumerate}
\item[(I)] $a<1,~ a+b\geq 1$ and $1+a+b-ab\leq c\leq 2$;
\item[(II)] $b\le 1$ and $2\le c.$
\end{enumerate}
\end{thm}



\section{Some lemmas}
In this section, we list some auxiliary lemmas
which play important roles in the proofs of
the main results.
The first result says that the class $\T$ is compact in the topology
of locally uniform convergence.

\begin{lem}\label{lem:compact}
Let $F_n~(n=1,2,\dots)$ be a sequence of functions in $\T.$
If $F_n$ converges to a function $F$ on $\Lambda$ pointwise.
Then $F\in\T.$
\end{lem}

\begin{pf}
By assumption, there is a sequence of probability measures $\mu_n$
on $[0,1]$ such that
\begin{equation}\label{eq:Fn}
F_n(z)=\int_0^1 \frac{d\mu_n(t)}{1-tz}, \quad z\in\Lambda.
\end{equation}
Since the set of (Borel) probability measures on $[0,1]$ is weakly compact,
there is a subsequence $\mu_{n_k}~(k=1,2,\dots)$ of $\mu_n$ which
weakly converges to a probability measure $\mu$ on $[0,1].$
Letting $k\to\infty$ in \eqref{eq:Fn} with $n=n_k,$ we obtain
a representation of $F$ in the form \eqref{eq:F} as required.
\end{pf}

As an immediate consequence, we obtain the following.

\begin{cor}
The class of universally convex functions is compact
in the topology of locally uniform convergence on $\Lambda.$
\end{cor}

In the proof, sometimes the hypergeometric function, which is involved,
may reduce to a polynomial so that the argument might break down.
In such a case, usually we can cover it by taking a suitable limit
of parameters, based on the above corollary.

We also need basic properties of functions in $\T$ in the following.

\begin{lem}\label{lem:Im}
Let $F\in\T.$ Then $\Im \big[zF(z)\big]\ge 0$ if $\Im z>0.$
\end{lem}

\begin{pf}
We first compute
$$
\Im\frac{z}{1-tz}=\Im \frac{z(1-t\bar z)}{|1-tz|^2}=\frac{\Im z}{|1-tz|^2}
$$
for $t\in[0,1].$
By integration with respect to the representing measure $d\mu(t)$ of $F$ over $[0,1],$
we obtain the conclusion.
\end{pf}

The next lemma is due to Ruscheweyh,
Salinas and Sugawa \cite{RusSalSu},
and is simplified by Liu and Pego \cite{LiuPego:2016}
(see also \cite[Lem.1.1]{LSW:2022}).

\begin{lem}\label{lem:RSS}
Let $F(z)$ be analytic on $\Lambda$. Then $F\in\T$
if and only if the following four conditions are fulfilled:
\begin{enumerate}
\item[(i)] $F(0)=1;$\medskip
\item[(ii)] $F(x)\in\R$ for $x\in(-\infty,1);$ \medskip
\item[(iii)] $\Im F(z)\ge 0$ for $\Im z>0;$ \medskip
\item[(iv)] $\displaystyle \limsup_{x\to+\infty}F(-x)\ge 0$.
\end{enumerate}
\end{lem}

\begin{rem}\label{rem:T}
By the integral form of functions in $\T$ we have indeed $F(x)$ is positive and
non-decreasing in $-\infty<x<1$ for every $F\in\T.$
In particular, the limit $\alpha=\lim_{x\to+\infty}F(-x)$ exists and satisfies
$0\le\alpha\le F(0)=1.$
Note also that $\alpha=1$ occurs only when $F(z)\equiv 1.$
We remark further that $G=(F-\alpha)/(1-\alpha)\in\T$ if $\alpha<1.$
Hence, $F=\alpha+(1-\alpha)G.$
Note that the class $\T$ is convex; namely, $aF+(1-a)G\in\T$ whenever
$F,G\in\T$ and $0\le a\le 1.$
In particular, $1-a+aF\in\T$ for $F\in\T$ and $0\le a\le 1$ because
the constant function $1$ belongs to $\T.$
\end{rem}

It is a good position to give an account on the proof of Theorem C.

\begin{pf}[Proof of Theorem C]
Put $g(z)=z\Gauss(a,2;c;z)$ and recall $f(z)=z\Gauss(a,1;c;z).$
Then $zf'(z)=g(z)$ (see \cite[p.~1368]{Kustner:2007}) and thus
$$
\Psi(z)=1+\frac{zf''(z)}{f'(z)}=\frac{zg'(z)}{g(z)}.
$$
Under the assumption, Theorem 1.8 in \cite{RusSalSu} implies
that $g$ is universally starlike; in other words, the function
$\Psi(z)=zg'(z)/g(z)$ belongs to the class $\T.$
Hence, we conclude that $1+zf''(z)/[2f'(z)]=(1+\Psi(z))/2$ belongs to $\T,$
too.
\end{pf}

The next lemma asserts that a nontrivial transformation preserves the class $\T.$

\begin{lem}[$\text{\cite[Lem.2.5]{Wang2021}}$]
\label{ratio}
Let $F\in\T.$ Then the function $\dfrac1{(1-z)F(z)}$ belongs to $\T.$
\end{lem}

The behavior of the hypergeometric function $\Gauss(a,b;c;z)$
at $z=1$ is completely different according to the
sign of $a+b-c$ as follows (see also \cite{AbramowitzStegun:1965}, \cite{OLBC:2010},
\cite{Ponn:1997}):

\begin{lem}\label{lem:1-}
Let $a,b,c\in\R$ with $c\ne 0,-1,-2,\dots.$
\begin{enumerate}
\item
If $c-a-b>0,$ the limit of $\Gauss(a,b;c;z)$ exists as $z\to 1$
in the domain $\Lambda=\C\setminus[1,+\infty)$ and it is given by
$$
\Gauss(a,b;c;1)=\frac{\Gamma(c)\Gamma(c-a-b)}{\Gamma(c-a)\Gamma(c-b)}.
$$
\item
If $c-a-b=0,$ as $z\to 1$ in $\Lambda$
$$
\Gauss(a,b;c;z)=\frac{\Gamma(c)}{\Gamma(a)\Gamma(b)}\log\frac1{1-z}+O(1).
$$
\item
If $c-a-b<0,$ as $z\to 1$ in $\Lambda$
$$
\Gauss(a,b;c;z)=\frac{\Gamma(c)\Gamma(a+b-c)}{\Gamma(a)\Gamma(b)}(1-z)^{c-a-b}
+O(|1-z|^{c-a-b+\varepsilon})
$$
for some number $\varepsilon>0.$
\end{enumerate}
\end{lem}

By using this, we can show the following fact.

\begin{lem}\label{lem:pres}
Let $0<a\le b\le c.$
Then
\begin{equation}\label{eq:pres1}
\lim_{x\to+\infty}\frac{\Gauss(a+1,b;c;-x)}{\Gauss(a,b;c;-x)}=0.
\end{equation}
If $b<c,$
\begin{equation}\label{eq:pres3}
\lim_{x\to+\infty}\frac{x\Gauss(a+1,b+1;c+1;-x)}{\Gauss(a+1,b;c;-x)}=+\infty.
\end{equation}
If furthermore $b\le a+1$ in addition to $b<c,$
\begin{equation}\label{eq:pres2}
\lim_{x\to+\infty}\frac{\Gauss(a+1,b+1;c+1;-x)}{\Gauss(a+1,b;c;-x)}=0.
\end{equation}
\end{lem}

\begin{pf}
We first recall the general formula (see \cite[15.3.4]{AbramowitzStegun:1965}):
$$
\Gauss(a,b;c;z)=(1-z)^{-a}\Gauss(a,c-b;c;\tfrac{z}{z-1})
$$
and the symmetry relation $\Gauss(a,b;c;z)=\Gauss(b,a;c;z).$
We now have
$$
\phi(x):=\frac{\Gauss(a+1,b;c;-x)}{\Gauss(a,b;c;-x)}
=\frac{(1+x)^{-a-1}\Gauss(a+1,c-b;c;\frac{x}{x+1})}%
{(1+x)^{-a}\Gauss(a,c-b;c;\frac{x}{x+1})}
$$
for $x>0.$
In what follows, we shall use the symbol $f(x)\sim g(x)$ to mean that $f(x)/g(x)$ tends
to a positive real number as $x\to +\infty.$
By Lemma \ref{lem:1-}, when $b<c,$ we have
$$
\phi(x)\sim\begin{cases}
x^{-1} & \quad (1<b-a), \\
x^{-1}\log x & \quad (b-a=1), \\
x^{a-b} & \quad (0<b-a<1), \\
1/\log x & \quad (b-a=0).
\end{cases}
$$
When $b=c,$ we have $\phi(x)=1/(1+x).$
In any case, $\phi(x)\to0$ as $x\to+\infty.$

Similarly, we compute
$$
\psi(x):=\frac{\Gauss(a+1,b+1;c+1;-x)}{\Gauss(a+1,b;c;-x)}
=\frac{(1+x)^{-a-1}\Gauss(a+1,c-b;c+1;\frac{x}{x+1})}%
{(1+x)^{-a-1}\Gauss(a+1,c-b;c;\frac{x}{x+1})}.
$$
By Lemma \ref{lem:1-} again, when $b<c,$ we have
$$
\psi(x)\sim\begin{cases}
1 & \quad (1<b-a), \\
1/\log x & \quad (b-a=1), \\
x^{b-a-1} & \quad (0<b-a<1), \\
x^{-1}\log x & \quad (b-a=0).
\end{cases}
$$
It is now easy to see that $x\psi(x)\to+\infty$ as $x\to+\infty,$
which means \eqref{eq:pres3}.
If $b\le a+1,$ the first case does not occur.
Thus we conclude that $\psi(x)\to 0$ as $x\to+\infty,$
which proves \eqref{eq:pres2}.
\end{pf}

The above \eqref{eq:pres1} yields the following formula.

\begin{lem}\label{lem:pres2}
Let $0<a\le b\le c$ and $f$ denote the shifted hypergeometric function
$z\Gauss(a,b;c;z).$
Then
$$
\lim_{x\to-\infty}\frac{x f''(x)}{f'(x)}=-a.
$$
\end{lem}

\begin{pf}
We denote by $h(z)$ the function $\Gauss(a+1,b;c;z)/\Gauss(a,b;c;z).$
In view of the formula (3.2) in \cite{Wang2021}, we have the formula
\begin{equation}\label{eq:W}
\frac{z f''(z)}{f'(z)}=\frac{2-c+(a+b-1)z}{1-z}
+\frac{c-2+(1-a)(1-b)z}{(1-z)(1-a+a h(z))}.
\end{equation}
Since $\phi(x)=h(-x),$ \eqref{eq:pres1} in Lemma \ref{lem:pres} now implies
$$
\lim_{x\to-\infty}\frac{x f''(x)}{f'(x)}
=-(a+b-1)-\frac{(1-a)(1-b)}{1-a}=-a.
$$
\end{pf}

In particular, we observe that for the shifted hypergeometric function
$f(z)=z\Gauss(a,b;c;z)$ with $0<a\le b\le c,$ this lemma
yields $\alpha=\lim_{x\to+\infty}F(-x)=1-a/2$ for $F(z)=1+zf''(z)/[2f'(z)].$
Thus the condition $a\le 2$ is necessary for $f$ to be universally convex
as we observed in Remark \ref{rem:T} (see also Theorem B).
We will indeed show that $\Psi(z)=1+zf''(z)/f'(z)$ belongs to $\T$
in the proof of Theorem \ref{thm:convex} so that the condition $a\le 1$ will be needed.

We also need some information about the asymptotic behaviour of the pre-Schwarzian derivative as
$x\to 1^-.$

\begin{lem}\label{lem:M2}
Let $0<a\le b\le c.$
Then
$$
\lim_{x\to1^-}\frac{\Gauss(a+1,b+1;c+1;x)}{\Gauss(a+1,b;c;x)}=\frac{c}{\max\{b,c-a-1\}}.
$$
\end{lem}

\begin{pf}
We put $H(x)=\Gauss(a+1,b+1;c+1;x)$ and $G(x)=\Gauss(a+1,b;c;x)$ and let
$\delta=(c+1)-(a+1)-(b+1)=c-(a+1)-b=c-a-b-1.$
When $\delta>0,$ by Lemma \ref{lem:1-} (1), we obtain
$$
\lim_{x\to1^-}\frac{H(x)}{G(x)}=\frac{\Gamma(c+1)\Gamma(\delta)}{\Gamma(c-a)\Gamma(c-b)}
\cdot\frac{\Gamma(c-a-1)\Gamma(c-b)}{\Gamma(c)\Gamma(\delta)}
=\frac{c}{c-a-1}.
$$
When $\delta=0,$ Lemma \ref{lem:1-} (2) implies
$$
\lim_{x\to1^-}\frac{H(x)}{G(x)}=\frac{\Gamma(a+b+2)}{\Gamma(a+1)\Gamma(b+1)}
\cdot\frac{\Gamma(a+1)\Gamma(b)}{\Gamma(a+b+1)}
=\frac{a+b+1}{b}
=\frac{c}{b}.
$$
When $\delta<0,$ Lemma \ref{lem:1-} (3) implies
$$
\lim_{x\to1^-}\frac{H(x)}{G(x)}=\frac{\Gamma(c+1)\Gamma(-\delta)}{\Gamma(a+1)\Gamma(b+1)}
\cdot\frac{\Gamma(a+1)\Gamma(b)}{\Gamma(c)\Gamma(-\delta)}
=\frac{c}{b}.
$$
Since $\delta\ge 0$ if and only if $b\le c-a-1,$ the above computations
are summarized to the assertion.
\end{pf}

The following lemma describes refined asymptotic behaviour of $\Gauss(a, b; c; z)$
when $c<a+b<c+1$.
The first assertion is essentially contained in \cite[Lemma 2.5]{PonVuor:1997}.

\begin{lem}\label{lem:asym}
Let $a,~b,~c$ be positive numbers with $a\le b$ and $\delta=a+b-c\in(0,1).$
Define a series $\{\sigma_n\}$ as the coefficients of the power series expansion of
\begin{equation}\label{positive}
R(z)=\Gauss(a, b; c; z)-A(1-z)^{-\delta}=\sum_{n=0}^\infty \sigma_n z^n,
\quad
A=\frac{\Gamma(c)\Gamma(\delta)}{\Gamma(a)\Gamma(b)}.
\end{equation}
Then $\sigma_n>0$ for all $n\ge0$ if $a<c<b,$ and $\sigma_n<0$ for all $n\ge0$
if either $c<a$ or $c>b.$
Also, the following relation is satisfied:
\begin{equation}\label{eq:B}
\sum_{n=0}^\infty \sigma_n=\frac{\Gamma(c)\Gamma(-\delta)}{\Gamma(c-a)\Gamma(c-b)}=:B.
\end{equation}
In particular, $|R(z)|\le |B|$ for $z\in\D$ and
$\Gauss(a, b; c; z)=A(1-z)^{-\delta}+O(1)$ as $z\to 1$ in $\D.$
\end{lem}

\begin{pf}
We compute $\sigma_n$ as follows:
$$
\sigma_n=\frac{(a)_n(b)_n}{(c)_n n!}-A\frac{(\delta)_n}{n!}
=\frac{(\delta)_n(Q_n-A)}{n!},
\quad\text{where}\quad
Q_n=\frac{(a)_n(b)_n}{(c)_n (\delta)_n}.
$$
First we note that
$$
\lim_{n\to\infty}Q_n= \frac{B(c,\delta)}{B(a,b)}=A,
$$
where $B(x,y)$ denotes the beta function
(this follows from Stirling's formula, see \cite[p.~283]{PonVuor:1997}).
Since
$$
\frac{Q_{n+1}}{Q_n}-1
=\frac{(a+n)(b+n)}{(c+n)(\delta+n)}-1
=\frac{(a-c)(b-c)}{(c+n)(\delta+n)},
$$
$Q_n$ is increasing and thus $Q_n<A$ when $(a-c)(b-c)>0$
and $Q_n$ is decreasing and thus $Q_n>A$ when $a<c<b.$
The first part has been proved.

To show the second part, we recall the formula
\begin{equation}\label{eq:TR}
\Gauss(a,b; c; z)=B\,\Gauss(a,b; \delta+1; 1-z)+A(1-z)^{-\delta}\Gauss(c-a,c-b; 1-\delta; 1-z)
\end{equation}
for $z\in\C\setminus((-\infty,0]\cup[1,+\infty))$
(see \cite[15.3.6]{AbramowitzStegun:1965}).
Therefore,
\begin{align*}
R(z)
&=B\,\Gauss(a,b; \delta+1; 1-z)+A(1-z)^{-\delta}\big(\Gauss(c-a,c-b; 1-\delta; 1-z)-1\big) \\
&=B\,\Gauss(a,b; \delta+1; 1-z)+O\big((1-z)^{1-\delta})
\end{align*}
and hence $R(z)\to B$ as $z\to1$ in $\D.$
Since the coefficients $\sigma_n$ have the same sign, it is easy to see that the series $\sum \sigma_n$
converges to $B,$ which completes the proof.
\end{pf}

\begin{lem}$($\cite[Lem.3.2]{SugawaWang:2016}$)$
\label{lem:conv}
Let $\Omega$ be an unbounded convex domain in $\C$ whose
boundary is parametrized positively by a Jordan curve $w(t)=u(t)+iv(t),~
0<t<1,$ with $w(0^+)=w(1^-)=\infty.$
Suppose that $u(0^+)=+\infty$ and that $v(t)$ has a finite limit
as $t\to0^+.$
Then $v(t)\le v(0^+)$ for $0<t<1.$
\end{lem}

The next two lemmas describe the
properties of the ratio of two hypergeometric functions
in different aspects.

\begin{lem}[$\text{\cite[Thm.1.5]{Kustner:2002},  \cite[ p.337-339 and Thm.69.2]{Wall:anal}}$]
\label{lem:Kust}
If $-1\leq a\leq c$ and $0\leq b\leq c\not=0$, then the three functions
$$
\frac{\Gauss(a+1,b;c;z)}{\Gauss(a,b;c;z)}, \quad
\frac{\Gauss(a+1,b+1;c+1;z)}{\Gauss(a,b;c;z)}
\aand
\frac{\Gauss(a+1,b+1;c+1;z)}{\Gauss(a+1,b;c;z)}
$$
belong to the class $\T.$
\end{lem}



\section{Proofs of the main results}

We will use the same symbols as in Lemma \ref{lem:asym} in the first two proofs.

\begin{pf}[Proof of Theorem \ref{thm:1}]
We begin with the proof of part (I).
Let $f(z)=z\Gauss(a,b;c;z),~ z\in\D.$
Then
$$
f(z)-A(1-z)^{-\delta}=z\big(R(z)+A(1-z)^{-\delta}\big)-A(1-z)^{-\delta}
=zR(z)-A(1-z)^{1-\delta}
$$
and therefore
$$
\big|f(z)-A(1-z)^{-\delta}\big|\le |zR(z)|+A|1-z|^{1-\delta}\le |B|+A\,2^{1-\delta}=\varepsilon.
$$
Since $A(1-z)^{-\delta}\in S(\delta)$ for $z\in\D,$ we obtain $f(z)\in S_\varepsilon(\delta).$

Next we show part (II).
Let $g(z)=z\Gauss(a,b;c;z^2),~ z\in\D.$
Then
$$
g(z)-A(1-z^2)^{-\delta}=zR(z^2)-A(1-z)^{1-\delta}(1+z)^{-\delta}.
$$
Hence, for $z\in\D$ with $\Re z\ge 0,$ we have $|1+z|\ge 1$ and
$$
\big|g(z)-A(1-z^2)^{-\delta}\big|\le |zR(z^2)|+A|1-z|^{1-\delta}|1+z|^{-\delta}
< |B|+A\,2^{1-\delta}=\varepsilon.
$$
Since $A(1-z^2)^{-\delta}\in S^*(\delta)$ for $z\in\D,$ we obtain $g(z)\in S_\varepsilon^*(\delta).$
Recall that the set $S_\varepsilon^*(\delta)$ is symmetric with respect to the origin.
Noting that $g$ is an odd function, we readily see that $g(z)\in S_\varepsilon^*(\delta)$
for $\Re z\le 0$ as well.
Hence $g(\D)\subset S_\varepsilon^*(\delta)$ follows.
\end{pf}

\begin{pf}[Proof of Theorem \ref{thm:sector}]
By Lemma \ref{lem:asym},
the function $\Gauss(a, b; c; z)$ is decomposed into
$\Gauss(a, b; c; z)=A(1-z)^{-\delta}+R(z).$
Let $\Omega$ be the image domain of $\D$
under the function $f(z)=z\Gauss(a, b; c; z)$.
Let $\gamma(\theta)=x(\theta)+iy(\theta)=f(e^{i\theta})$,
$\gamma_1(\theta)=x_1(\theta)+iy_1(\theta)=Ae^{i\theta}(1-e^{i\theta})^{-\delta}$ and
$\gamma_2(\theta)=x_2(\theta)+iy_2(\theta)=e^{i\theta}R(e^{i\theta})$
for $0< \theta<2\pi$.
Note that the curves $\gamma(\theta)$, $\gamma_1(\theta)$
and $\gamma_2(\theta)$ are symmetric in the real axis.
Therefore, we may restrict our attention to the case when $0<\theta\le\pi.$
We first determine the asymptotic line of the Jordan arc $\gamma(\theta)$ as $\theta\to 0^{+}$.
An elementary computation yields that
$$
x_1(\theta)=\frac{A\cos \left(\theta+\frac{\delta}{2}(\pi-\theta)\right)}%
{\left(2 \sin \frac{\theta}{2}\right)^\delta}
\aand
y_1(\theta)=\frac{A\sin \left(\theta+\frac{\delta}{2}(\pi-\theta)\right)}%
{\left(2 \sin \frac{\theta}{2}\right)^\delta}.
$$
Thus $\displaystyle\lim_{\theta\to 0^{+}}x_1(\theta)=\displaystyle\lim_{\theta\to 0^{+}}y_1(\theta)=+\infty$
because of $0<\delta<1$.
On the other hand, the curve $x_2(\theta)+iy_2(\theta)$ is bounded and, as $\theta\to 0^{+},$
$x_2(\theta)\to R(1)=B$ and $y_2(\theta)\to 0.$
Therefore
\begin{eqnarray*}
\lim_{\theta\to 0^{+}}\frac{y(\theta)}{x(\theta)}
=\lim_{\theta\to 0^{+}}\frac{y_1(\theta)+y_2(\theta)}
{x_1(\theta)+x_2(\theta)}
=\tan\left(\frac{\delta\pi}{2}\right)
\end{eqnarray*}
and
\begin{align*}
& \, \quad  \lim_{\theta\to 0^{+}}\left[y(\theta)-\tan\left(\frac{\delta\pi}{2}\right) x(\theta)\right]\\
&=\lim_{\theta\to 0^{+}}\left[y_1(\theta)-\tan\left(\frac{\delta\pi}{2}\right) x_1(\theta)\right]+
\lim_{\theta\to 0^{+}}\left[y_2(\theta)-\tan\left(\frac{\delta\pi}{2}\right) x_2(\theta)\right]\\
&=\lim_{\theta\to 0^{+}}\frac{A\sin \left(\theta-\frac{\delta}{2}\theta\right)}
{\cos\left(\frac{\delta\pi}{2}\right) \left(2 \sin \frac{\theta}{2}\right)^\delta}
-B\tan\left(\frac{\delta\pi}{2}\right)
=-B\tan\left(\frac{\delta\pi}{2}\right).
 \\
\end{align*}
For brevity, we write $\phi=\delta\pi/2.$
We have obtained the asymptotic line $l: y=(\tan\phi) (x-B)$ of the Jordan curve $\gamma(\theta)$
as $\theta\to 0^{+}.$
By making a rotation $w=e^{-i\phi}(z-B)+B$ about the point $B$,
the convex domain $\Omega$ and the asymptotic line $l$ are transformed
to a convex domain $\Omega'$ and the line $y=0$.
If we write $w(\theta)=u(\theta)+iv(\theta)=e^{-i\phi}(\gamma(\theta)-B)+B,$
we have $u(0^+)=+\infty$ and $v(\theta)=-(\sin\phi) (x(\theta)-B)+(\cos\phi) y(\theta)
=(\cos\phi)[y(\theta)-(\tan\phi) (x(\theta)-B)]\to 0$ as $\theta\to 0^+.$
We now apply Lemma \ref{lem:conv} to see that
the convex domain $\Omega'$ lies below the line $y=0$,
which implies that the convex domain $\Omega$
lies below the asymptotic line $l$.
By symmetry, we now conclude that
the convex domain $\Omega$ is contained in the sector $|\arg (z-B)|<\phi=\delta\pi/2$
as required.
\end{pf}

We should look at the pre-Schwarzian derivative $f''/f'$ of $f(z)=z\Gauss(a,b;c;z)$
in order to obtain universal convexity of $f$ on $\Lambda$.
Based on the expression \eqref{eq:W} in terms of $\Gauss(a+1,b;c;z)/\Gauss(a,b;c;z),$
we will derive two other representations in terms of
 $\Gauss(a+1,b+1;c+1;z)/\Gauss(a,b;c;z)$ and
$\Gauss(a+1,b+1;c+1;z)/\Gauss(a+1,b;c;z)$ in Propositions \ref{prop:integral} and
\ref{prop:integral2}, respectively.

\begin{prop}\label{prop:integral}
Let $a,~b$ and $c$ be positive numbers
with $a\leq b\leq c$ and $a\le 1$ and put $f(z)=z\Gauss(a,b;c;z).$
If the inequality
\begin{equation}\label{eq:ab}
(2-c)(c+ab-a-b-1)=(2-c)(1-a)(1-b)-(2-c)^2\ge 0
\end{equation}
holds, then
$$
\frac{f''(z)}{f'(z)}=
\frac{a+b-1}{1-z}+
\left(1-a-b+\frac{2ab}{c}\right)\Phi_1(z),
$$
for a function $\Phi_1\in\T.$
\end{prop}

\begin{pf}
We first consider the special case when $(1-a)(1-b)=0.$
Note that \eqref{eq:ab} then implies $c=2.$
For instance, if $a=1,$ then we can easily see that $f'(z)=(1-z)^{-b}$
and $f''(z)/f'(z)=b/(1-z).$
Thus the conclusion is confirmed in this case.
Therefore, in the following, we will assume that either $(1-a)(1-b)\ne0.$
Since $a\le 1,$ we have thus $a<1.$

We set $\tau=ab/c~(>0)$ for the sake of brevity and make preliminary observations.
When $c+ab-a-b-1\ge0$ and $2-c\ge0,$ we have
$$
1-a-b+2\tau\ge 1-a-b+2(1+a+b-c)/c=\frac{(1+a+b)(2-c)}{c}\ge 0.
$$
When $c+ab-a-b-1\le0$ and $2-c\le0,$ we have similarly
$1-a-b+2\tau\le 0.$
Therefore, the assumption \eqref{eq:ab} implies that
$c+ab-a-b-1, 1-a-b+2\tau, 1-b, (1-a)(1-b)-(2-c)$
and $2-c$ have the same sign (admitting to be 0).

We next write $F(z)=\Gauss(a,b;c;z)$, $G(z)=\Gauss(a+1,b;c;z)$
and $H(z)=\Gauss(a+1,b+1;c+1;z)$ for convenience.
We infer from the contiguous relation
\begin{equation}\label{contiguous}
G(z)-F(z)=\frac{b}{c}\,zH(z)
\end{equation}
that \eqref{eq:W} can now be rearranged as
\begin{align*}
\frac{z f''(z)}{f'(z)}
&
=\frac{2-c+(a+b-1)z}{1-z}
+\frac{c-2+(1-a)(1-b)z}{(1-z)(1+\tau zH(z)/F(z))}\\
&=\frac{(a+b-1)z}{1-z}+
\frac{z}{1-z}
      \frac{(1-a)(1-b)+(2-c)\tau H(z)/F(z)}{1+\tau zH(z)/F(z)}.
\end{align*}
Hence,
\begin{equation*}
\frac{ f''(z)}{f'(z)}
=\frac{a+b-1}{1-z}+\frac{1-a-b+2\tau}{(1-z)M(z)},
\end{equation*}
where
\begin{equation}\label{eq:M(z)}
M(z)=\frac{(1-a-b+2\tau)[1+\tau z w(z)]}{(1-a)(1-b)+(2-c)\tau w(z)}
\aand w(z)=\frac{H(z)}{F(z)}.
\end{equation}
We show now that $M\in\T,$ which yields that $\Phi_1(z)=1/[(1-z)M(z)]$ belongs to $\T$
by Lemma \ref{ratio}.
For this,  it is sufficient to verify conditions (i)-(iv) in Lemma \ref{lem:RSS}.
Condition (i): $M(0)=1$ is easy to check.
Since $0<a\leq b\leq c$, we can apply Lemma \ref{lem:Kust} to obtain $w\in\T;$
that is, $w$ is represented  by
$$
w(z)=\int_0^1\frac{d\mu(t)}{1-tz},\quad z\in\Lambda,
$$
for a probability measure $\mu$ on $[0,1].$
We now check condition (ii).
Since $(1-a)(1-b)$ and $2-c$ have the same sign 
and $w(x)>0$ for $x\in(-\infty,1),$
the denominator $(1-a)(1-b)+(2-c)\tau w(x)$ of $M(x)$ is either positive for all $x<1$
or negative for all $x<1.$
Hence, the denominator of $M(x)$ does not vanish and thus
$M(x)\in\R$ for $x<1,$ which confirms condition (ii).
In order to verify condition (iii), we compute
the imaginary part of $M(z)$ as follows:
\begin{eqnarray*}
&&\Im M(z)\left|(1-a)(1-b)+(2-c)\tau w(z)\right|^2\\
&=&(1-a-b+2\tau)
\Im\left[\left[1+\tau z w(z)\right]
\left((1-a)(1-b)+(2-c)\tau\overline{w(z)}\right)\right]\\
&=&(1-a-b+2\tau)\Im \left((1-a)(1-b)\tau z w(z)-(2-c)\tau w(z)
+(2-c)\tau^2z\left|w(z)\right|^2\right)\\
&=&(1-a-b+2\tau)\tau\Im \left(\int_{0}^{1}\frac{(1-a)(1-b) z-(2-c)}{1-tz}d\mu(t)
+(2-c)\tau z\left|w(z)\right|^2\right)\\
&=&(1-a-b+2\tau)\tau y\left(\int_{0}^{1}\frac{(1-a)(1-b)-(2-c)t}{|1-tz|^2}d\mu(t)
+(2-c)\tau \left|w(z)\right|^2\right),
\end{eqnarray*}
for $z=x+iy$.
Making use of the preliminary observations above, we see that
$(1-a-b+2\tau)[(1-a)(1-b)-(2-c)t]\geq 0$ for $0\le t\le 1$ and
$(1-a-b+2\tau)(2-c)\geq 0,$
by which we verify condition (iii).
By \eqref{eq:pres3} and the contiguous relation \eqref{contiguous},
we have
\begin{equation*}
\lim_{x\to+\infty}xw(-x)
=\lim_{x\to+\infty}\frac{xH(-x)}{G(-x)+bxH(-x)/c}
=\lim_{x\to+\infty}\frac{cxH(-x)/G(-x)}{c+bxH(-x)/G(-x)}
=\frac{c}{b},
\end{equation*}
and, in particular, $w(-x)\to 0$ as $x\to +\infty.$
Thus, when $a,b\ne1,$ condition (iv) is checked by computing
\begin{align*}
\lim_{x\to+\infty}M(-x)
&=(1-a-b+2\tau)\lim_{x\to+\infty}
\frac{1-\tau xw(-x)}
{(1-a)(1-b)+(2-c)\tau w(-x)}\\
&= \frac{(1-a-b+2\tau)(1-a)}{(1-a)(1-b)}=\frac{1-a-b+2\tau}{1-b}
\geq0.
\end{align*}
Hence, we have now concluded that $M(z)$ belongs to the class
$\T,$ as required.
\end{pf}

\begin{prop}\label{prop:integral2}
First let $a,~b$ and $c$ be positive numbers
with $a\le b\le 1$ and $c\ge 2$.
Then, the pre-Schwarzian derivative of
the shifted hypergeomegtric function $f(z)=z\Gauss(a,b;c;z)$
is expressed by
$$
\frac{zf''(z)}{f'(z)}=a\Phi_2(z)-a
$$
for a function $\Phi_2\in\T.$
\end{prop}

\begin{pf}
We use the same notation as in the previous proof.
In addition, set $\sigma=(a-1)b/c$
and $\sigma'=(c-a-1)b/c$.
For a while, we assume that $0<a\le b<1<2<c.$
By the proof of Theorem 1.2 of \cite{Wang2022} we have
\begin{align*}
\frac{zf''(z)}{f'(z)}&=-a+
\frac{a(c-a-1)}{(1-a)(1-z)}+
\frac{a}{1-a}\cdot\frac{2-c+(1-a)(b-1)z}{(1-z)(1+\sigma z H(z)/G(z))} \\
&=-a+\frac{a}{1-z}\frac{1+(b-1)z-\sigma'zH(z)/G(z)}{1+\sigma zH(z)/G(z)}\\
&=-a+\frac{a}{(1-z)M(z)},
\end{align*}
where $M=M_1/M_2$ and
$$
M_1(z)=1+\sigma z w_1(z), \quad
M_2(z)=1+(b-1)z-\sigma' zw_1(z),  \quad
w_1(z)=\frac{H(z)}{G(z)}.
$$
We show now that $M\in\T,$
which implies that $\Phi_2(z)=1/[(1-z)M(z)]$ belongs to $\T$
by virtue of Lemma \ref{ratio}.
(Note that the proof of Theorem 1.2 of \cite{Wang2022} was flawed \cite{Wang2022Cor}
but the above formula is valid.)
Lemma \ref{lem:Kust} implies that
$$
w_1(z)=\int_0^1\frac{d\mu(t)}{1-tz},\quad z\in\Lambda,
$$
for a probability measure $\mu$ on $[0,1].$
In order to show that $M\in\T,$ we only need to check conditions (i)-(iv)
in Lemma \ref{lem:RSS}.
Condition (i) follows by definition.
We show now (ii).
The point is to see $M_2(x)=1+(b-1)x-\sigma'x w_1(x)\ne0$ for $x<1.$
Since
$$
xw_1(x)=\int_0^1\frac{x}{1-tx}d\mu(t)
$$
is non-decreasing in $-\infty<x<1$, $b-1<0$ and $\sigma'>0$,
the function $M_2(x)$ is decreasing in $x<1.$
we need only to show that $M_2(1^-)\ge 0.$
By Lemma \ref{lem:M2},
\begin{align*}
\lim_{x\to1^-}M_2(x)&=b-\sigma'\lim_{x\to 1^-}\frac{H(x)}{G(x)} \\
&=b-\frac{(c-a-1)b}{c}\cdot\frac{c}{\max\{b,c-a-1\}} \\
&=\max\{1+a+b-c,0\}\ge 0
\end{align*}
and condition (ii) has been verified.

Next we show condition (iii).
As before, by making use of Lemma \ref{lem:Im},
we have
\begin{align*}
&\quad \ \Im M(z)\left|1+(b-1)z-\sigma' zw_1(z)\right|^2\\
&=
\Im\left[\left[1+\sigma z w_1(z)\right]
\left(1+(b-1)\overline{z}-\sigma' \overline{zw_1(z)}\right)\right]\\
&=(1-b)y+(\sigma+\sigma')\Im[zw_1(z)]+(b-1)\sigma|z|^2\Im w_1(z)\\
&> 0
\end{align*}
for $z=x+iy$ with $y>0$, since all the coefficients
$1-b,~\sigma+\sigma'=(c-2)b/c,~(b-1)\sigma$ are positive.

Finally, we check (iv) in Lemma \ref{lem:RSS}.
By \eqref{eq:pres2}, we have
$$
M(-x)=\frac{M_1(-x)}{M_2(-x)}
=\frac{1/x-\sigma w_1(-x)}{1/x+(1-b)+\sigma'w_1(-x)}
\to0
$$
as $x\to+\infty$
and thus (iv) is confirmed.
Thus the assertion has been proved when $a\le b<1$ and $2<c.$
The general case is verified by using, for example, the approximation
$a_n=a-1/n, b_n=b-1/n, c_n=c+1/n ~(n$: large enough) with
the compactness property of the class
$\T$ (see Lemma \ref{lem:compact}).
Now the proof is complete.
\end{pf}

\begin{pf}[Proof of Theorem \ref{thm:convex}]
We will employ the same symbols as in the proof of Proposition \ref{prop:integral}.
Thus $f(z)$ means the shifted hypergeometric function $z\Gauss(a,b;c;z)$.
By Theorem B, it is enough to
prove that the function $1+zf''(z)/[2f'(z)]=(1+\Psi(z))/2$ belongs to the class $\T,$
where $\Psi(z)=1+zf''(z)/f'(z).$
By convexity of $\T,$ it suffices to show that
$\Psi$ belongs to $\T.$
Thus we only need to verify conditions (i)--(iv) in Lemma \ref{lem:RSS} for $\Psi.$
At this point, we note that (iv) is guaranteed by Lemma \ref{lem:pres2},
which says that $\Psi(-x)\to 1-a\ge 0$ as $x\to+\infty.$

\noindent
\textbf{Case (I)}:
Note that the assumptions of Proposition \ref{prop:integral} are satisfied in this case.
Thus the function $\Psi$ is expressed as
$$
\Psi(z)=1+
\frac{(a+b-1)z}{1-z}+(1-a-b+2\tau) z\Phi_1(z), \quad z\in\Lambda.
$$
Conditions (i) and (ii) are clearly satisfied.
Condition (iii) is easily verified by using Lemma \ref{lem:Im}.


\noindent
\textbf{Case (II)}:
We make use of Proposition \ref{prop:integral2}.
All the assumptions of the proposition are satisfied in this case.
Recall the formula there:
$$
\Psi(z)=1+\frac{zf''(z)}{f'(z)}=1-a+a\Phi_2(z)
$$
for some $\Phi_2\in\T.$
Hence $\Psi\in\T$ by convexity of $\T.$
\end{pf}



\def\cprime{$'$} \def\cprime{$'$} \def\cprime{$'$}
\providecommand{\bysame}{\leavevmode\hbox
to3em{\hrulefill}\thinspace}
\providecommand{\MR}{\relax\ifhmode\unskip\space\fi MR }
\providecommand{\MRhref}[2]{%
  \href{http://www.ams.org/mathscinet-getitem?mr=#1}{#2}
} \providecommand{\href}[2]{#2}

\end{document}